\documentclass[11pt,reqno]{amsart}
\usepackage[dvips]{epsfig}
\RequirePackage{amsmath}
\usepackage{epsfig}
\usepackage{graphicx,epsfig}
\usepackage{amscd,amsthm,amsfonts,amsopn,amssymb,verbatim}
\usepackage{color}
\usepackage{times}
\numberwithin{equation}{section}
\theoremstyle{remark}

\def\be{\begin{equation}}
\def\ee{\end{equation}}
\def\vp{\varphi}

\allowdisplaybreaks
\def\be{\begin{equation}}
\def\ee{\end{equation}}
\def\vp{\varphi}

\allowdisplaybreaks
\begin{document}
\title
[]
{The essential centre of the mod a diagonalization ideal commutant of an
$n$-tuple of commuting hermitian operators}
\author
{J.~Bourgain and D.V.~Voiculescu}
\address
{School of Mathematics, Institute for Advanced Study, Princeton, NJ 08540}
\address
{Department of Mathematics, UC Berkeley, Berkeley, CA 94720-3840}
\email
{bourgain@math.ias.edu}
\thanks{Research supported in part by NSF Grants DMS 1301619 for first author and DMS 1001881 for second author.}

\abstract
We show that for a commuting $n$-tuple of hermitian operators, with
perfect spectrum, the essential centre of the algebra of operators
commuting with the $n$-tuple mod a diagonalization ideal arises from the
$C^*$-algebra of the $n$-tuple.
This answers a question for normal operators and the Hilbert-Schmidt
class connected to $K$-theory for almost normal operators.
\endabstract
\date{\today}

\maketitle
Let $\tau = (T_1, \ldots, T_n)$ be an $n$-tuple of commuting hermitian
operators on a complex separable infinite-dimensional Hilbert space
$\mathcal H$, let $(J, | \ \ |_J)$ be a normed ideal of compact operators
on $\mathcal H$ in which the finite rank operators $\mathcal R(\mathcal
H)$ are dense and let $\mathcal E(\tau, J)$ be the algebra of bounded
operators $X\in \mathcal B(\mathcal H)$ so that $[X, T_j]\in J, 1\leq
j\leq n$.
Let further $\mathcal K(\tau; J)=\mathcal E(\tau; J)\cap \mathcal
K(\mathcal H)$, where $\mathcal K(\mathcal H)$ denotes the compact
operators and let $\mathcal E/\mathcal K(\tau; J)=\mathcal E(\tau;
J)/\mathcal K(\tau; J)$.

Roughly, the main result of this note is about a situation when the
cenre of $\mathcal E/\mathcal K(\tau; J)$ consists of the image of
$C^*(\tau)$, the $C^*$-algebra of $\tau$, in $\mathcal E/\mathcal
K(\tau; J)$.
In case $n=2$ and $J=\mathcal C_2$, the Hilbert-Schmidt class, this
answers a question in 6.2 of \cite{3} and also provides a generalization.
The motivation for the question about the essential centre came from the
$K$-theory problems studied in \cite{3}.

The main assumption will be that $\tau$ can be diagonalized modulo $J$,
that is that there is a hermitian $n$-tuple $\delta =(D_1, \ldots, D_n)$
which is diagonal in some orthonormal basis of $\mathcal H$, so that
$D_j-T_j\in J$, $1\leq j\leq n$.
By the results of \cite {2}, this is equivalent to the requirement that
$$
k_J(\tau)=\underset{{A\in \mathcal R^+_1 (\mathcal H)}\atop
{A\uparrow I}} {\lim \inf} \ \max_{1\leq j\leq n}|[A, T_j]|_J=0
$$
where the lim inf is w.r.t. the natural order on the set of finite rank
positive contractions $\mathcal R_1^+(\mathcal H)$ on $\mathcal H$.

Also by \cite {2}, if $k_J(\tau)=0$ we can choose $\delta$ so that we
have the equality of spectra $\sigma(\tau)=\sigma(\delta)$ and that
$|T_j -D_j|_J<\varepsilon$, $1\leq j\leq n$.

In what follows $\mathcal B/\mathcal K(\mathcal H)=\mathcal B(\mathcal
H)/\mathcal K(\mathcal H)$ is the Calkin algebra and \hfill\break
$p:\mathcal B(\mathcal H)\to \mathcal B/\mathcal K(\mathcal H)$ denotes the
canonical homomorphism.

Endowed with the norm
$$
||| X ||| =\Vert X\Vert +\max_{1\leq j\leq n} |[X, T_j]|_J
$$
the algebra $\mathcal E(\tau, J)$ is easily seen to be an involutive
Barnach algebra with isometric involution.

\medskip
\noindent
{\bf Theorem.}
{\sl Let $\tau =(T_1, \ldots, T_n)$ be an $n$-tuple of commuting hermitian
operators on $\mathcal H$ such that $k_J(\tau)=0$.

a) The algebraic isomorphism of $*$-algebras $\mathcal E/\mathcal
K(\tau; J) \sim p\big(\mathcal E(\tau; J)\big)$ is an isometric
isomorphism of Banach algebras with involution.
In particular $\mathcal E/\mathcal K(\tau, J)$ is a $C^*$-algebra.

b) Assume the spectrum $\sigma(\tau)$ is a perfect set.
Then, under the isomorphism $\mathcal E/\mathcal K(\tau; J)\sim
p\big(\mathcal E(\tau; J)\big)$, the centre $\mathcal Z\big(\mathcal E/
\mathcal K(\tau; J)\big)$ corresponds to $p\big (C^*(\tau)\big)$.
}

\begin{proof}
a) The proof is along the same lines as the proof of Proposition 5.3 in
\cite {3}.
Since $k_J(\tau)=0$, there are $A_m\in\mathcal R_1^+(\mathcal H)$ so that
$A_m\uparrow I$ and $|[A_m, T_j]|_J\to 0$, $1\leq j\leq n$ as $m\to
\infty$.
If $X\in \mathcal E(\tau, J)$ we have 
$$
\begin{aligned}
&\underset{m\to\infty} {lim \, sup} |[(I-A_m)X, T_j]|_J\leq\\
&\leq \underset{m\to\infty} {lim \, sup} \Vert X\Vert \ |[I-A_m, T_j]|_J+\\
&\underset{m\to\infty} {lim\, sup} |(I-A_m)[X, T_j]|_J=0.
\end{aligned}
$$
Here, the last limsup equals zero since $[X, T_j]\in J$ and we assumed
that the finite rank operators $\mathcal R(\mathcal H)$ are dense in
$J$.

It follows that in $\mathcal E/\mathcal K(\tau; J)$ we have
$$
\begin{aligned}
&||| X+\mathcal K(\tau; J)||| \leq\\
&\leq \underset{m\to\infty}{\text lim sup} |||(I-A_m)X|||\leq\\
&\leq \underset{m\to\infty}{\text lim sup}\Vert (I-A_m)X\Vert+\\
& \underset{m\to\infty}{\text lim sup} \max_{1\leq  j\leq n} |[(I-A_m) X, T_j]|_J=\\
&=\Vert p(X)\Vert.
\end{aligned}
$$
The opposite inequality is obvious, which concludes the proof of a).
\medskip

b) Denoting commutants by primes, we have $\mathcal E(\tau; J)\supset(\tau)'$.

Hence if $X+\mathcal K(\tau; J) \in \mathcal Z\big(\mathcal E/\mathcal K (\tau; J)\big)$ we have
$[X, (\tau)']\subset\mathcal K(\mathcal H)$ and by a result of \cite {1}, this implies
$X\in(\tau)'' +\mathcal K(\mathcal H)$.
Since $(\tau)''\in\mathcal E(\tau; J)$ we infer $X\in (\tau)'' +\mathcal K(\tau; J)$.

Since $\sigma(\tau)$ is a perfect set we have $\sigma(\tau) =\sigma\big(p(\tau)\big)$.
Also, note that if $\tilde\tau = (\tilde T_1, \ldots, \tilde T_n)$ is another $n$-tuple of
commuting hermitian operators such that \break
$T_j -\tilde T_j\in J, 1\leq j\leq n$ we have $\mathcal
E(\tilde \tau, J) =\mathcal E(\tau; J)$ and we may replace $\tau$ by $\tilde \tau$.
Thus, using \cite{2} we may assume $\tau$ is diagonal in some orthonormal  basis and
$\sigma(\tau)=\sigma \big(p(\tau)\big)$ continues to be satisfied.
Further, it suffices to deal with the case of $X=X^*$ and assume $X=\vp(T_1, \ldots, T_n)$
where $\vp:\sigma(\tau)\to \mathbb R$ is a bounded Borel function.

Let $\Delta \subset\sigma(\tau)$ be the dense subset of $\sigma(\tau)$ which is the joint pure
point spectrum and let $d$ be the distance on $\sigma(\tau)$ corresponding to the $\ell^\infty$-norm
on $\ell^\infty(\{1, \ldots, n\})\sim \mathbb R^n$.

To prove part b) of the theorem, we must show that $\vp$ can be chosen to be continuous and this in
turn is easily seen to be equivalent to showing that $\vp(k_m)$ is convergent as $m\to\infty$ 
whenever $(k_m)_{m\in\mathbb N}\subset \Delta$ is a Cauchy sequence.
Further, passing to subsequences, it is sufficient to show this in case $1> d(k_m, k_{m+1})> 10 d(k_{m+1}, k_{m+2})$.
Again passing to subsequences and replacing $\vp$ by some $a\vp+b$, we conclude that under the previous 
assumptions we must show that 
$$
|\vp(k_{2r})|\leq 1/10, |\vp(k_{2r+1})-1|\leq 1/10
$$
for all $r\in\mathbb N$, implies we can find $Y\in \mathcal E(\tau; J)$ so that 
$$
[Y, \vp(\tau)]\not\in \mathcal K(\mathcal H).
$$

For each $m\in \mathbb N$ let $e_m$ be an eigenvector of $\tau$ in the eigenspace for the $n$-tuple of eigenvalues
$k_m$, so that $(e_m)_{m\in \mathbb N}$ is an orthonormal system.

We define
$$
Y=\sum^\infty_{m=1} \langle \cdot, e_m\rangle e_{m+1}.
$$
Then $[Y, T_j]$ has $s$-numbers majorized by $d(k_m, k_{m+1})$, $m\in\mathbb N$ and since these are $\leq C10^{-m}$
we have $[Y, T_j]\in\mathcal C_1\subset J$, so that $Y\in\mathcal E(\tau; J)$.

On the other hand
$$
[Y,\vp(\tau)]=\sum^\infty_{m=1} \big(\vp(k_m)-\vp(k_{m+1})\big) \langle\cdot, e_m\rangle e_{m+1}
$$
which is not compact since $\vp(k_m)-\vp(k_{m+1})$ does not converge to zero as $m\to\infty$
\end{proof}

\end{document}